\documentclass[12pt]{amsart}
\usepackage[utf8]{inputenc}
\usepackage{amsmath}
\usepackage{xcolor}
\usepackage{amssymb}
\newtheorem{theorem}{Theorem}
\newtheorem*{theorem*}{Theorem}

\newtheorem{corollary}{Corollary}[section]
\newtheorem{lemma}{Lemma}

\newtheorem*{acknowledgements*}{Acknowledgements}

\makeatletter
\def\blfootnote{\gdef\@thefnmark{}\@footnotetext}
\makeatother

\def\house#1{\setbox1=\hbox{$\,#1\,$}%
\dimen1=\ht1 \advance\dimen1 by 2pt \dimen2=\dp1 \advance\dimen2 by 2pt
\setbox1=\hbox{\vrule height\dimen1 depth\dimen2\box1\vrule}%
\setbox1=\vbox{\hrule\box1}%
\advance\dimen1 by .4pt \ht1=\dimen1
\advance\dimen2 by .4pt \dp1=\dimen2 \box1\relax}

\begin{document}
\title[Quadratic non-residues and  non-primitive roots satisfying a coprimality condition]{Quadratic non-residues and  non-primitive roots satisfying a coprimality condition}
\author{Jaitra Chattopadhyay, Bidisha Roy, Subha Sarkar and R. Thangadurai}
   \address[Jaitra Chattopadhyay, Bidisha Roy, Subha Sarkar and R. Thangadurai]{Harish-Chandra Research Institute, HBNI\\
Chhatnag Road, Jhunsi\\
Allahabad 211019, India}
\email[Jaitra Chattopadhyay]{jaitrachattopadhyay@hri.res.in}
\email[Bidisha Roy]{bidisharoy@hri.res.in}
\email[Subha Sarkar]{subhasarkar@hri.res.in}
\email[R. Thangadurai]{thanga@hri.res.in}

\begin{abstract}
Let $q\geq 1$ be any integer and let $ \epsilon \in [\frac{1}{11}, \frac{1}{2})$ be a given real number. 
In this short note, we prove that for all primes $p$ satisfying 
$$
p\equiv 1\pmod{q}, \quad \log\log p > \frac{\log 6.83}{\frac{1}{2}-\epsilon} \mbox{ and } \frac{\phi(p-1)}{p-1} \leq \frac{1}{2} - \epsilon,
$$ 
there exists a quadratic non-residue $g$ which is not a primitive root modulo $p$  such that $gcd\left(g, \frac{p-1}{q}\right) = 1$. 
\end{abstract}


\maketitle

\section{Introduction}

Let $p$ be an odd prime number. We know that there are exactly $\frac{p-1}{2}$ quadratic residues as well as non-residues modulo $p$. It is a well known fact that the multiplicative group $\left(\mathbb{Z}/p\mathbb{Z}\right)^*$ is cyclic (see \cite{apostol}). An element is called a {\it primitive root} modulo $p$ if it is a generator of this cyclic group.

The distribution of quadratic residues, non-residues and primitive roots is a very fundamental area in number theory and has been a topic of immense interest to mathematicians for centuries.  In 2010,  Levin, Pomerance and Soundararajan \cite{sound-discrete} proved the following result.

\begin{theorem}\label{sound}
For all  prime numbers $p\geq 5$, there exists a primitive root $g$ modulo $p$ which satisfies the condition $gcd(g, p-1) =1$.
\end{theorem}

Levin, Pomerance and Soundararajan \cite{sound-discrete} considered this problem to tackle a particular case of an important problem in computational   number theory, namely, discrete log problem. More precisely, they prove Theorem \ref{sound} to tackle the fixed point discrete log problem which states that {\it for a given primitive root $g$ in $(\mathbb{Z}/p\mathbb{Z})^*$, does there exists an integer $t\in [1, p-1]$ such that $g^t \equiv t\pmod{p}$?} Indeed, Theorem \ref{sound} solves the fixed point discrete log problem affirmatively. 

In this article, we deal with the similar  problem  for quadratic non-residues which are not primitive roots (for further reference on this related problem, see  \cite{tha1}, \cite{jarso} and \cite{tha3}). For notational convenience, we abbreviate  `a quadratic non-residue which is not a primitive root' by QNRNP modulo $p$.  More precisely, we prove the following result.

\begin{theorem}\label{thm1}
Let $q \geq 1$ be an integer and $\epsilon \in [\frac{1}{11},\frac{1}{2})$. Let $p$ be a prime satisfying 
$$
p\equiv 1\pmod{q}, \quad \log\log p > \frac{\log 6.83}{\frac{1}{2}-\epsilon} \mbox{ and } \frac{\phi(p-1)}{p-1}\leq\frac{1}{2}-\epsilon
.$$ 
Then there exists an integer $g$ satisfying $1 < g < p-1$ and $gcd\left(g,\frac{p-1}{q}\right)=1$ such that $g$ is a \text{QNRNP} modulo $p$. In particular, when $q =1$, there exists an integer $g$ with $1<g<p-1$ and $gcd(g, p-1)=1$ such that $g$ is a \text{QNRNP} modulo $p$.
 \end{theorem}
 
 In the statement of Theorem \ref{thm1}, one of the conditions on $p$ is  a natural condition. If  $\displaystyle\frac{\phi(p-1)}{p-1} = \frac{1}2$, then one can easily check that every non-residue modulo $p$ is a primitive root modulo $p$. The condition $\displaystyle \frac{\phi(p-1)}{p-1} \leq \frac{1}{2} -\epsilon$ makes sure that $p-1$ has enough odd prime factors and hence abundance of QNRNP residues modulo $p$. 
  As an application, we solve the fixed point discrete log problem for the cyclic subgroup generated by a QNRNP as follows. 
 
 \begin{corollary}\label{cor1}
 Let $\epsilon \in [\frac{1}{11},\frac{1}{2})$ be a real number. Let $p$ be a prime satisfying 
$$
\log\log p > \frac{\log 6.83}{\frac{1}{2}-\epsilon} \mbox{ and } \frac{\phi(p-1)}{p-1}\leq \frac{1}{2}-\epsilon
.$$ 
Then  there is a QNRNP $g$ and an integer $x \in [1,p-1]$ such that $x$  is  QNRNP and 
$g^x \equiv x \pmod{p}$.
\end{corollary}

In \cite{sound-discrete}, first they proved their result for all large primes and used the computations to check their result for small primes. However, computations may be cumbersome in our result stated in  Theorem \ref{thm1} because of various parameters. 
  
\section{Preliminaries}

Let $\mu_{p-1}$ stand for the multiplicative group of $(p-1)$-th roots of unity. Let $g \in \{1, \ldots, p-1\}$ be a primitive root modulo $p$ and let $\chi : (\mathbb{Z}/p\mathbb{Z})^{*} \rightarrow \mu_{p-1}$ be a character modulo $p$ such that $\chi$ is a generator of the dual group of $(\mathbb{Z}/p\mathbb{Z})^*$. For all integers $\ell$ with $0\leq \ell \leq p-2$, we denote $\chi_\ell = \chi^\ell$ a character modulo $p$ and $\chi_0$ is the principal character.

Suppose $\chi (g)=\eta$. Since $\chi$ is a generator of the dual group of $(\mathbb{Z}/p\mathbb{Z})^*$ and $g$ is a primitive root modulo $p$, we get that $\eta$ is a primitive $(p-1)$-th root of unity. 

Following \cite{tha1}, we define
$$
\beta_{\ell}(p-1)=\displaystyle \sum_{\substack{1 \leq i \leq p-1 \\i  \ \ \ {\textit odd} \\(i,p-1) > 1}}(\eta^{i})^{\ell} \qquad \mbox{ and } \qquad  \alpha_{\ell}(p-1) = \sum_{\substack{1 \leq i \leq p-1 \\ (i,p-1) =1}}(\eta^{i})^{\ell},
$$
where $\alpha_{\ell}(p-1)$ is known as {\it Ramanujan sums}.  

Now, we list some basic lemmas and results which will be useful to us in course of the proof of Theorem \ref{thm1}.

\begin{lemma}\label{lem01}\cite{tha1}
For  all integers $\ell$ with $0<\ell<p-1$, we have 
$$
\beta_{\ell}(p-1)=-\alpha_{\ell}(p-1).
$$
\end{lemma}

\begin{lemma}(characteristic function for $QNRNP$'s)\label{lem02} \cite{tha1}
For any $x\in (\mathbb{Z}/p\mathbb{Z})^*$, we have 
$$
\displaystyle \sum_{\ell=0}^{p-2}\beta_{\ell}(p-1)\chi_{\ell}(x) =\left\{
\begin{array}{ll}
p-1; & \mbox{ if }  \mbox{x is a QNRNP},\\
0; & \mbox{ otherwise}.
\end{array}
\right.$$
\end{lemma}

Now,  we shall state some basic results as follows.

\begin{lemma} \label{lem1}
\noindent \begin{enumerate}
 
\item Let $\omega(n)$ denote the number of distinct prime divisors of $n$. Then
we have
$$
\omega(p-1) \leq (1.385) \frac{\log p}{\log\log p}$$
for all primes $p\geq 5$. (See for instance, \cite{sandor}).

\item For any positive integer $n$, let $\mu(n)$ denote the M\"{o}bius function. Then, we have 
$$
\sum_{d|n}\mu(d) = \left\{\begin{array}{ll}
1;&\mbox{ if } n = 1,\\
0;&\mbox{ if } n > 1.
\end{array}\right.
$$
(See for instance, \cite{apostol}).

\item For any odd prime $p$ and any divisor $q$ of $p-1$, we have
$$
\sum_{d|\frac{p-1}{q}} |\mu(d)| = 2^{\omega\left(\frac{p-1}q\right)}.
$$
(See for instance, \cite{sza1}). 
\end{enumerate}
\end{lemma}

\begin{lemma}\label{lem4} \cite{sza1}
We have, 
$$
\displaystyle \sum_{\ell=1}^{p-2}|\alpha_{\ell}(p-1)|=2^{\omega(p-1)}\phi (p-1).
$$
\end{lemma}

The following result is a standard theorem to estimate a character sum over an interval which can be found in \cite{apostol}. 

\begin{theorem}\label{polya1} (P\'{o}lya-Vinogradov) 
Let $p$ be any odd prime and $\chi$ be a non-principal character modulo $p$. Then, for any integers $0\leq M < N\leq p-1$, we have,  
$$
\left|\sum_{m=M}^N \chi(m)\right| \leq \sqrt{p} \log p.$$
\end{theorem}

\section{Proof of Theorem \ref{thm1}}

Let $q\geq 1$ be a given integer and let $ \epsilon \in [\frac{1}{11}, \frac{1}{2})$ be also given. Now, we consider all primes $p\equiv 1\pmod{q}$ with $\frac{\phi(p-1)}{p-1}\leq \frac{1}{2}-\epsilon$. By Dirichlet's prime number theorem, we can see that there are infinitely many  such primes. 

By Lemma \ref{lem02},  for any integer $m$, we let, 
$$
f(m):=\frac{1}{p-1}\displaystyle \sum_{\ell=0}^{p-2}\beta_{\ell}(p-1)\chi_{\ell}(m)=\left\{\begin{array}{ll}
1; & \mbox{ if } m  \mbox{ is a QNRNP},\\
0; & \mbox{ for otherwise}.
\end{array}\right.
$$
By letting $N_p := \displaystyle \sum_{\substack{m=1 \\ (m, \frac{p-1}{q})=1}}^{p-1}f(m)$, we see that   $N_p$ counts the number of QNRNP's in $\{1,\ldots,p-1\}$ which are relatively prime with $\frac{p-1}{q}$. To finish the proof of Theorem 1, it suffices to prove that $N_{p} \geq 1$ for all primes $p > \exp \exp \frac{\log 6.83}{\frac{1}{2}-\epsilon}$ satisfying $p\equiv 1\pmod{q}$ and $\frac{\phi(p-1)}{p-1}\leq \frac{1}{2}-\epsilon$.  Therefore, we consider  
\begin{eqnarray*}
N_p = \displaystyle \sum_{\substack{m=1 \\ \left(m,\frac{p-1}{q}\right)=1}}^{p-1}f(m)&=&\frac{1}{p-1}\displaystyle \sum_{\substack{m=1 \\ \left(m,\frac{p-1}{q}\right)=1}}^{p-1}\sum_{\ell=0}^{p-2}\beta_{\ell}(p-1)\chi_{\ell}(m)\\  
&=& \frac{1}{p-1}\sum_{\ell=0}^{p-2}\beta_{\ell}(p-1)\displaystyle \sum_{\substack{m=1 \\ \left(m,\frac{p-1}{q}\right)=1}}^{p-1}\chi_{\ell}(m)\\ 
&=& \frac{1}{p-1}\left(\beta_{0}(p-1)q\phi\left(\frac{p-1}{q}\right)+\sum_{\ell=1}^{p-2}\beta_{\ell}(p-1)\displaystyle \sum_{\substack{m=1 \\ \left(m,\frac{p-1}{q}\right)=1}}^{p-1}\chi_{\ell}(m)\right),
\end{eqnarray*}
where we have used the fact that the number of integers $m$ in $\{1, \ldots, p-1\}$ such that $\left(m, \frac{p-1}{q}\right)= 1$ is $q\phi\left(\frac{p-1}{q}\right)$. Let us define 
 \begin{eqnarray*}
E_{p} &:=& N_{p}-\frac{1}{p-1}\beta_{0}(p-1)q\phi\left(\frac{p-1}{q}\right) = \frac{1}{p-1}\sum_{\ell=1}^{p-2}\beta_{\ell}(p-1)\displaystyle \sum_{\substack{m=1 \\ \left(m,\frac{p-1}{q}\right)=1}}^{p-1}\chi_{\ell}(m). 
\end{eqnarray*}
In order to prove $N_p\geq 1$, we need to get an upper bound for $E_{p}$. For that, we need to estimate 
$\displaystyle \sum_{\substack{m=1 \\ \left(m,\frac{p-1}{q}\right)=1}}^{p-1}\chi_{\ell}(m)$ and $\displaystyle \frac{1}{p-1}\sum_{\ell=1}^{p-2}\beta_{\ell}(p-1)$ separately. First we consider the sum  $\displaystyle \sum_{\substack{m=1 \\ \left(m,\frac{p-1}{q}\right)=1}}^{p-1}\chi_{\ell}(m)$ as follows. For a given integer $\ell$ with $1\leq \ell\leq p-2$, we have 
\begin{eqnarray*}
\displaystyle \sum_{\substack{m=1 \\ \left(m,\frac{p-1}{q}\right)=1}}^{p-1}\chi_{\ell}(m) &=& \sum_{m=1}^{p-1}\chi_{\ell}(m)\sum_{d \mid \left(m,\frac{p-1}{q}\right)}\mu(d) = \sum_{d \mid \frac{p-1}{q}}\mu(d)\displaystyle \sum_{\substack{t=1}}^{\frac{p-1}{d}}\chi_{\ell}(d)\chi_{\ell}(t)\\ 
&=& \sum_{d \mid \frac{p-1}{q}}\mu(d)\chi_\ell(d)\sum_{t=1}^{\frac{p-1}{d}}\chi_{\ell}(t),
\end{eqnarray*}
by Lemma \ref{lem1} (2).
Hence, by Theorem \ref{polya1} and Lemma \ref{lem1} (3), we get 
\begin{eqnarray*}
\left|\sum_{\substack{m=1 \\ \left(m,\frac{p-1}{q}\right)=1}}^{p-1}\chi_{\ell}(m)\right| &\leq & \sum_{d \mid \frac{p-1}{q}}\left|\mu(d)\right| \left|\sum_{t=1}^{\frac{p-1}{d}}\chi_{\ell}(t)\right|  \leq  2^{\omega\left(\frac{p-1}{q}\right)}\sqrt{p}\log p.
\end{eqnarray*}
Also, by Lemma \ref{lem01} and Lemma \ref{lem4}, we see that 
\begin{eqnarray*}
\left|\sum_{\ell=1}^{p-2}\beta_{\ell}(p-1)\right| & \leq & \sum_{\ell=1}^{p-2}\left|\beta_{\ell}(p-1)\right|                                    
                                             \leq  \sum_{\ell=0}^{p-2}\left|\alpha_{\ell}(p-1)\right|
                                             = 2^{\omega (p-1)}\phi(p-1).
\end{eqnarray*}
Thus, using the above two estimates, we get,
\begin{eqnarray} \label{epestimate}
\left| E_{p}\right| &=& \left| N_{p}-\frac{1}{p-1}\beta_{0}(p-1)q\phi\left(\frac{p-1}{q}\right)\right| \leq \frac{1}{p-1} \sum_{\ell=1}^{p-2}|\beta_{\ell}(p-1)|\cdot \left|\displaystyle \sum_{\substack{m=1 \\ \left(m,\frac{p-1}{q}\right)=1}}^{p-1}\chi_{\ell}(m)\right| \nonumber \\
 & \leq & 2^{\omega \left(\frac{p-1}{q}\right) + \omega (p-1)} \frac{\phi (p-1)}{p-1} \sqrt{p} \log p.
\end{eqnarray}
Observe that \eqref{epestimate}  implies
$$
-2^{\omega \left(\frac{p-1}{q}\right) +\omega (p-1)}\frac{\phi (p-1)}{p-1}\sqrt{p}\log p \leq N_{p}-\frac{q\phi \left(\frac{p-1}{q}\right)}{(p-1)}\beta_{0}(p-1),
$$
which is equivalent to 
$$
N_{p} \geq \frac{\phi \left(\frac{p-1}{q}\right)}{\frac{p-1}{q}}\beta_{0}(p-1)- 2^{\omega \left(\frac{p-1}{q}\right)+\omega(p-1)}\frac{\phi (p-1)}{p-1}\sqrt{p}\log p.
$$
Thus to establish $N_{p}>0$, it is enough to show that,
$$
\frac{\phi \left(\frac{p-1}{q}\right)}{\frac{p-1}{q}}\beta_{0}(p-1)-2^{\omega \left(\frac{p-1}{q}\right)+\omega (p-1)}\frac{\phi (p-1)}{p-1}\sqrt{p}\log p>0,
$$

which is equivalent to showing that 
\begin{equation}\label{eq2}
\beta_{0}(p-1)>2^{\omega \left(\frac{p-1}{q}\right)+\omega (p-1)} \frac{\phi(p-1)}{q\phi \left(\frac{p-1}{q}\right)}\sqrt{p}\log p.
\end{equation}

Now, it is clear that 
\begin{equation}\label{eq200}
\phi(p-1) \leq q \phi\left(\frac{p-1}{q}\right) \ \iff \ \ \frac{\phi(p-1)}{q\phi \left(\frac{p-1}{q}\right)} \leq 1.
\end{equation}

Since $\omega\left(\frac{p-1}{q}\right)\leq \omega(p-1)$, by \eqref{eq2} and \eqref{eq200}, it is enough to prove that 
\begin{equation}\label{eq201}
\beta_0(p-1) > 4^{\omega(p-1)} \sqrt{p}\log p,
\end{equation}
for primes $p > \exp \exp \frac{\log 6\cdot83}{\frac{1}{2}-\epsilon}$ satisfying $\frac{\phi (p-1)}{p-1} \leq \frac{1}{2} - \epsilon.$ 

Let $p$ be a prime satisfying $p> \exp \exp \frac{\log 6.83}{\frac{1}{2}-\epsilon}$.  Therefore, we get
\begin{equation}\label{eq3}
p^{\frac{1}{2}-\epsilon}>p^{\frac{\log 6.83}{\log \log p}}.
\end{equation}

By Lemma \ref{lem1} (1), we also know that 
$$
\omega (p-1) \leq 1.385 \frac{\log p}{\log \log p}.
$$

Therefore, we get 
$$4^{\omega(p-1)}\leq 4^{1.385 \frac{\log p}{\log \log p}}\leq 6.83^{\frac{\log p}{\log \log p}}=p^{\frac{\log 6.83}{\log \log p}}.
$$

Hence, from \eqref{eq3}, we have, 
$$p^{\frac{1}{2}-\epsilon}>4^{\omega(p-1)} \iff p^{1-\epsilon}(\log p)>4^{\omega (p-1)}\sqrt{p}\log p.
$$

In order to prove (\ref{eq201}), it is enough to show that 
\begin{equation}\label{eq4}
\beta_{0}(p-1)>p^{1-\epsilon}\log p
\end{equation}

for all primes $p > \exp \exp \frac{\log 6.83}{\frac{1}{2}-\epsilon}$ satisfying $\frac{\phi(p-1)}{p-1} \leq \frac{1}{2}-\epsilon$. 

\smallskip
 
Note that the condition 
$$
\frac{\phi(p-1)}{p-1} \leq \frac{1}{2}-\epsilon  \ \iff \ \epsilon(p-1) \leq \frac{p-1}{2} - \phi(p-1)  = \beta_0(p-1).
$$ 

Therefore, to prove \eqref{eq4}, it is enough to prove that  $\epsilon(p-1)\geq p^{1-\epsilon}\log p$ for all primes $\displaystyle p > \exp \exp \frac{\log 6.83}{\frac{1}{2}-\epsilon}$. 

\smallskip

Since $\epsilon \in [\frac{1}{11}, \frac{1}{2})$, we write $\epsilon = \frac{1}{c}$ for some real number $c$ with $2 < c\leq 11$ and note that 
$$
\log \log p > \frac{\log 6.83}{\frac{1}{2} -\epsilon} > 3.84 \times 1.22 > 4.68 \qquad \mbox{and} \qquad \log p > e^{4.68} > 107.7.
$$
In order to prove $\epsilon(p-1) \geq p^{1-\epsilon}\log p$ for all primes $\displaystyle p > \exp \exp \frac{\log 6.83}{\frac{1}{2}-\epsilon}$, it is enough to prove that 
$$
\frac{p}{1.1} > \frac{1}{\epsilon} p^{1-\epsilon}\log p \iff p > (1.1c)^c (\log p)^c \iff \log p > c\log(1.1c) +c \log \log p.
$$
Since we know $e^x/x \geq 22$ for all $x \geq 4.68$, we apply for $x =\log \log p$ and see that\\ $\log p > 2c \log\log p$ for all $2 < c \leq 11$. Hence, it is enough to prove that $$\log p > c\log(1.1) + c\log c +\frac{\log p}{2} \iff  \log p > 2c \log(1.1) + 2c\log c.$$ 

 Since $c \leq 11$, we see that 
$$2c\log(1.1) + 2c\log c \leq 22\log(1.1)+ 22\log 11 \leq 54.86 < 107.7< \log p.
$$
Thus the inequality in  \eqref{eq4} holds true,   which completes the proof of the theorem. $\hfill\Box$

\section{Proof of Corollary \ref{cor1}}

By Theorem \ref{thm1}, there is a QNRNP $x$ modulo $p$ satisfying $x \in [1,p-1]$ and $gcd(x,p-1) = 1$. Let $y$ be the
multiplicative inverse of $x$ modulo $p-1$.  Put $g = x^y$. Then note that $g$ is also a QNRNP modulo $p$. Hence, we get 
$g ^x \equiv x^{xy} \equiv x \pmod{p}$.

\bigskip

\noindent{\bf Acknowledgement} We are thankful to the referee for pointing out a lacuna in the previous version  and for suggesting important references. Also, we are grateful to Professor John Loxton for a careful reading which made us clear some ambiguity in the proof.

\end{document}